\documentclass[11pt,a4paper]{article}

\usepackage{amsmath}
\usepackage{amssymb}
\usepackage{amsthm}
\usepackage{amsfonts}
\usepackage{mathtools}
\usepackage{mathrsfs}

\usepackage[a4paper,margin=1in]{geometry}

\usepackage{microtype}
\usepackage{xurl}

\usepackage{amsmath}
\usepackage{mathtools}

\allowdisplaybreaks

\usepackage[
    colorlinks=true,
    linkcolor=blue,
    citecolor=blue,
    urlcolor=blue
]{hyperref}

\usepackage{cite}
\usepackage[colorlinks=true,
            citecolor=blue,
            linkcolor=blue,
            urlcolor=blue]{hyperref}

\newtheorem{theorem}{Theorem}[section]

\newtheorem{proposition}[theorem]{Proposition}
\newtheorem{corollary}[theorem]{Corollary}

\theoremstyle{definition}
\newtheorem{definition}[theorem]{Definition}
\newtheorem{example}[theorem]{Example}

\theoremstyle{remark}
\newtheorem{remark}[theorem]{Remark}
\usepackage{mathabx}

\title{\textbf{{Generalized Cram\'er-Rao Inequalities Associated with\\ $(h,\tau)$-Divergence and $(h,\tau)$-Exponential Family}}}

\author{
Satyajit Dhadumia\\
\small \textit{Department of Mathematics, IIT Palakkad}\\
\small \texttt{Corresponding author:satya22dhadumia@gmail.com}
}

\date{}

\begin{document}

\maketitle

\begin{abstract}{

The Kullback–Leibler (KL) divergence and exponential family have laid the foundation for information geometry and have been widely generalized. Using Eguchi’s theory, the Fisher information metric can be derived from the KL divergence. Based on this metric, Amari and Nagaoka established an information-geometric proof of the classical Cramér–Rao (CR) inequality for finite-support models. Among the various generalizations of the KL divergence, this paper focuses on the $(h,\tau)$-divergence and establishes the associated CR inequality, referred to as the $(h,\tau)$-CR inequality, using information geometry. We also provide a proof of this inequality based on the Cauchy–Schwarz inequality, thereby extending its applicability to models with more general support sets. Under suitable regularity conditions, equality holds in the classical CR inequality if and only if the underlying family is exponential. However, for the $(h,\tau)$-CR inequality, this property does not generally extend to the $(h,\tau)$-exponential family. We therefore establish an attainable CR inequality for this family. This new inequality extends the attainment result of the classical CR inequality by the sufficient statistic of the exponential family to the $(h,\tau)$-exponential family.}

\medskip


\end{abstract}
\section{Introduction}
Let $(\mathcal{X},\mu)$ be a measurable space, and let $\mathcal{P}$ denote the set of all probability densities with respect to $\mu$. Throughout, $\mu$ is the counting measure when $\mathcal{X}$ is finite or countably infinite and $\mu$ is the n-dimensional Lebesgue measure when $\mathcal{X}\subseteq \mathbb{R}^n$ is a measurable set with positive Lebesgue measure.

Let $S \subseteq \mathcal{P}$ be a family of probability distributions parameterized by $\theta=[\theta_1,\cdots,\theta_k]^{\top}$, where $k\geq1$, such that
\begin{equation}\label{eqn:statistical model}
    S=\{p_{\theta}:\theta=[\theta_1,\cdots,\theta_k]^{\top}\in \Theta\},
\end{equation}
where $\Theta$ is an open subset of $\mathbb{R}^k$, called the parameter space. Let us also suppose that $\operatorname{Supp}(p_{\theta}) = \mathcal{X}$ for all $\theta\in\Theta$ and the mapping $\theta\mapsto p_{\theta}$ is injective. For example, consider the $n$-dimensional multivariate normal family
\begin{equation}
\label{eqn:multivariate_normal}
p_{\mu,\sigma^2}(\mathrm{x})
=
\frac{1}{(2\pi\sigma^2)^{n/2}}
\exp\left\{
-\frac{1}{2\sigma^2}
\sum_{i=1}^n(x_i-\mu)^2
\right\},\quad\mathrm{x}=[x_1,\ldots,x_n]^{\top}\in\mathbb{R}^n,
\end{equation}
where $\mu\in\mathbb{R}$ and
$\sigma^2>0$. In information geometry, $S$ can be viewed as a $k$-dimensional differentiable manifold whose points represent probability distributions. We call it as a statistical manifold. The family of n-dimensional multivariate normal family forms a 2-dimensional statistical manifold. 

The Cramér-Rao (CR) inequality states that, under suitable regularity conditions, if $\hat{r}(\mathrm{X})$ is an unbiased estimator of $r(\theta)=[r_1(\theta),\ldots,r_k(\theta)]^{\top}$, a $k$-dimensional differentiable function of $\theta$, then
\begin{equation}
    \label{eqn:CRLB}
    V_{\theta}[\hat{r}(\mathrm{X})]\succeq J(\theta)[G_n(\theta)]^{-1}J(\theta)^{\top},
\end{equation}
where $G_n(\theta)\coloneqq
\mathbb{E}_{\theta}[u_{\theta}(\mathrm{X})
u_{\theta}(\mathrm{X})^{\top}]$ is the \emph{Fisher information (FI) matrix}, with $\mathrm{X}=(X_1,\ldots,X_n)\sim p_{\theta}$, $u_{\theta}(\mathrm{x})={\partial}/{\partial{\theta}}\log p_{\theta}(\mathrm{x})$ is the score function  and $J(\theta)$ is the $k \times k$ Jacobian matrix of $r(\theta)$ \cite[eq.~(3.5)]{crrao1945}, \cite[Th.~7.3.9]{casella2024statistical}\\
(Note: Here and throughout, for two matrices $A_1$ and $A_2$, we write $A_1\succeq A_2$ to indicate that the matrix $A_1-A_2$ is positive semi-definite). The proof of \eqref{eqn:CRLB} is a direct application of the Cauchy–Schwarz inequality.

  If the FI-matrix $G_n(\theta)=[g_{ij}(\theta)]_{i,j=1}^k$ is a positive definite symmetric matrix for all $\theta$, then it defines a unique Riemannian metric (inner product):
  \begin{equation}
      g(\theta)\coloneqq\langle X,Y   \rangle=\sum_{i,j=1}^k g_{ij}(\theta)X^iY^j,
  \end{equation} for any two tangent vectors $X,\,Y\in T_{p_{\boldsymbol{\theta}}}(S)$ where $X=\sum_{i=1}^k X^i\partial_i,\,Y=\sum_{j=1}^k Y^{j}\partial_j$ with $X^i,\,Y^j\in \mathbb{R}$, and $\partial_i\coloneqq\frac{\partial}{\partial\theta_i}$, thereby transforming $S$ into a Riemannian manifold. This metric is known as the \emph{Fisher information metric}. This metric can also be obtained from the Kullback–Leibler (KL) divergence using Eguchi’s theory \cite{Eguchi1992}:
\begin{equation}
\label{eqn:Fisher_from_KL}
g_{ij}(\theta)
=
-\left.\partial_i\partial'_j
D_{KL}(p_\theta,p_{\theta'})
\right|_{\theta'=\theta},
\end{equation}
where $\partial'_j\coloneqq\frac{\partial}{\partial{\theta'}_j}$, and the KL divergence is defined as
\begin{equation}\label{eqn:KL-div}
D_{KL}(p_{\theta},p_{\theta'})
\coloneqq
\int_{\mathcal{X}} p_{\theta}(\mathrm{x})\log\frac{p_{\theta}(\mathrm{x})}{p_{\theta'}(\mathrm{x})}\,d\mu(\mathrm{x}).
\end{equation} Based on this metric, Amari-Nagaoka provided an information geometric proof of the inequality \eqref{eqn:CRLB} for models with finite support (see Section 2.5 of \cite{AmariHiroshi2000}).\\
Observe that, with $h(t)=t\log t$, $\tau(t)=t$, \eqref{eqn:KL-div} is a special case of $(h,\tau)$-divergence \cite[eq. (1.2)]{MATSUZOE2027131088}:
\begin{equation}
    \label{eqn:(h,tau)-divergence}
D_{(h,\tau)}(p_{\theta},p_{\theta'})\coloneqq\int \left\{h(\tau(p_\theta(\mathrm{x})))-h(\tau(p_{\theta'}(\mathrm{x})))-h'(\tau(p_{\theta'}(\mathrm{x})))[\tau(p_\theta(\mathrm{x}))-\tau(p_{\theta'}(\mathrm{x}))]\right\}d\mu(\mathrm{x}),
\end{equation}where $\tau$ is a smooth function on an open interval $I$ contained in $(0,\infty)$ and $h$ is a smooth function on $\tau(I)$ such that $\tau'>0$ on $I$ and $h''>0$ on $\tau(I)$. For $h(t)=\varphi(t)$ , $\tau(t)=t$, \eqref{eqn:(h,tau)-divergence} coincide with the Bregman divergence \cite{Roy2019DensityPower}:
\begin{equation*}
    D_{(\varphi)}(p_{\theta},p_{\theta'})\coloneqq\int \Big[{\varphi(p_\theta(\mathrm{x}))-\varphi(p_{\theta'}(\mathrm{x}))-\varphi'(p_{\theta'}(\mathrm{x}))\big(p_\theta(\mathrm{x})-p_{\theta'}(\mathrm{x})\big)\Big]}d\mu(\mathrm{x}),
\end{equation*} where $\varphi$ is strictly convex twice differentiable function on $(0,\infty)$.

In this paper, we establish the CR inequality associated with the $(h,\tau)$-divergence. To establish this inequality, we first apply Eguchi's theory to derive the corresponding Riemannian metric. We then employ Amari–Nagaoka's approach to obtain the desired inequality, referred to as the
$(h,\tau)$-CR inequality. Since Amari–Nagaoka's approach is restricted to finite-support models, the established inequality is applicable only to such models. Consequently, we also prove it using the Cauchy–Schwarz inequality, thereby extending its applicability to models with more general supports. Under suitable regularity conditions, equality holds in the classical CR inequality if and only if the underlying family is exponential \cite{fend1959attainment,wijsman1973attainment,muller1989attainment}. However, this attainment property does not generally extend to the $(h,\tau)$-CR inequality under an $(h,\tau)$-exponential family. We therefore establish an attainable CR inequality for this family. This new inequality extends the attainment result of the classical CR inequality by the sufficient statistic of the exponential family to the $(h,\tau)$-exponential family. The rest of the paper is organized as follows. In Section~\ref{sec:2}, we establish the $(h,\tau)$-CR inequality, from the $(h,\tau)$-divergence using an information-geometric approach. In Section~\ref{sec:3},
we extend this inequality to models with more general
support sets. In Section~\ref{sec:4}, we establish
a CR inequality that is attained under the
$(h,\tau)$-exponential family. Finally, we conclude
the paper with a summary in Section~\ref{sec:5}.
\section{A Generalized CR Inequality from $(h,\tau)$-Divergence using Information Geometry}\label{sec:2}
In this section, we assume that $\mathcal{X}$ is a finite set. Then, the probability simplex $\mathcal{P}$ (strictly positive probability distributions) forms a $|\mathcal{X}|-1$ dimensional statistical manifold. Let $ \theta = [\theta_1, \dots, \theta_k]^{\top}$ be a coordinate system of $\mathcal{P}$. By using Eguchi's theory, we can define a Riemannian metric $\widebar{g}^{(h,\tau)}(\theta)$ derived from the $(h,\tau)$-divergence $D_{(h,\tau)}$, characterized by the matrix $\widebar{G}_n^{(h,\tau)}(\theta)=[\widebar{g}_{ij}^{(h,\tau)}(\theta)]_{i,j=1}^k$ with
\begin{eqnarray}
\label{ij entry of phi fisher matrix}
\widebar{g}^{(h,\tau)}_{ij}(\theta) &\coloneqq& -D_{(h,\tau)}(\partial_i|\partial_j)\nonumber\\
    &\coloneqq& -\partial_i\partial'_j D_{(h,\tau)}(p_{\theta},p_{\theta'})\Big|_{\theta^{'}=\theta}\nonumber\\
     &=&\partial_i\partial_j'\Big\{\sum_{\mathrm{x}\in\mathcal{X}}h'(\tau(p_{\theta'}(\mathrm{x})))[\tau(p_\theta(\mathrm{x}))-\tau(p_{\theta'}(\mathrm{x}))]\Big\}\Big|_{\theta'=\theta}\nonumber\\
     &=&\sum_{\mathrm{x}}h''(\tau(p_{\theta}(\mathrm{x})))[\tau'(p_{\theta}(\mathrm{x}))]^2\frac{\partial p_{\theta}(\mathrm{x})}{\partial\theta}\frac{\partial p_{\theta}(\mathrm{x})}{\partial \theta}\nonumber\\
     &=&\sum_{\mathrm{y}}[h''(\tau(p_{\theta}(\mathrm{y})))[\tau'(p_{\theta}(\mathrm{y}))]^2]^{-1}\nonumber\\
     &&\mathbb{E}_{(h,\tau,\theta)}\left[\frac{1}{p_{(h,\tau,\theta)}(\mathrm{X})}\frac{\partial p_{\theta}(\mathrm{X})}{\partial\theta}\frac{1}{p_{(h,\tau,\theta)}(\mathrm{X})}\frac{\partial p_{\theta}(\mathrm{X})}{\partial\theta}\right],\nonumber
\end{eqnarray}where the expectation $\mathbb{E}_{{(h,\tau,\theta)}}$ is with respect to the $(h,\tau)$-escort distribution $$p_{(h,\tau,\theta)}(\mathrm{x})\coloneqq\frac{[h''(\tau(p_{\theta}(\mathrm{x})))[\tau'(p_{\theta}(\mathrm{x}))]^2]^{-1}}{\sum_{\mathrm{y}}[h''(\tau(p_{\theta}(\mathrm{y})))[\tau'(p_{\theta}(\mathrm{y}))]^2]^{-1}}\qquad {\mathrm{x}}\in\mathcal{X}.$$
By following the same setting as in \cite[Sec.~2.5]{AmariHiroshi2000}, we have 
\begin{equation}
    T_{p}^{(e)}(\mathcal{P})\coloneqq\left\{X^{(e)}\mid X\in T_{p}(\mathcal{P})\right\}=\left\{A\in \mathbb{R}^{\mathcal{X}}\mid \mathbb{E}_{p}{[A]}=0\right\},
\end{equation}where 
$X^{(e)}(\mathrm{x})\coloneqq X^{(m)}(\mathrm{x})/p(\mathrm{x}) \text{ for } \mathrm{x}\in \mathcal{X}$ with $X^{(m)}\in\mathcal{A}_0
\coloneqq\{A\in \mathbb{R}^{\mathcal{X}}\mid \sum_\mathrm{x} A(\mathrm{x})=0\}.$
For $p\in \mathcal{P}$, define
\begin{equation}
      X^{(h,\tau)}(\mathrm{x})=X^{(m)}(\mathrm{x})/p_{(h,\tau)}(\mathrm{x}),
\end{equation}where the escort $$p_{(h,\tau)}(\mathrm{x})\coloneqq\frac{[h''(\tau(p(\mathrm{x})))[\tau'(p(\mathrm{x}))]^2]^{-1}}{\sum_{\mathrm{y}}[h''(\tau(p(\mathrm{y})))[\tau'(p(\mathrm{y}))]^2]^{-1}}.$$
 $X^{(h,\tau)}$ is the $(h,\tau)$-representation of $X\in T_p(\mathcal{P})$. Let \begin{equation}
    T_p^{(h,\tau)}(\mathcal{P})\coloneqq\big\{X^{(h,\tau)}\mid X\in T_p(\mathcal{P})\big\}.\nonumber
\end{equation}
The following proposition gives another characterization for $T_p^{(h,\tau)}(\mathcal{P})$.

\medskip
\begin{proposition}
    \begin{equation}
        T_{p}^{(h,\tau)}(\mathcal{P}) = \left\{A \in \mathbb{R}^{\mathcal{X}}\mid \mathbb{E}_{p_{(h,\tau)}}[A] = 0\right\}.  
    \end{equation}
\end{proposition}
\begin{proof}\ Since $X^{(h,\tau)}(\mathrm{x})=X^{(m)}(\mathrm{x})/p_{(h,\tau)}(\mathrm{x})$. Hence $\mathbb{E}_{p_{(h,\tau)}}[X^{(h,\tau)}]=\sum_{x}X^{(m)}(\mathrm{x})=0$. On the other hand, if $A\in\mathbb{R}^{\mathcal{X}}$ such that $\mathbb{E}_{p_{(h,\tau)}}[A]=0$ then $\sum_\mathrm{x} A(\mathrm{x})p_{(h,\tau)}(\mathrm{x})=0$. This implies $A(\mathrm{x})p_{(h,\tau)}(\mathrm{x})=X^{(m)}(\mathrm{x})$. Thus, $A(\mathrm{x})=X^{(m)}(\mathrm{x})/p_{(h,\tau)}(\mathrm{x})\in T_p^{(h,\tau)}$.
\end{proof}
Since multiplying a Riemannian metric by a smooth, positive function yields another Riemannian metric, for any $X,Y\in T_{p}(\mathcal{P})$, we define the $(h,\tau)$-Fisher information metric at a point $p$ on $\mathcal{P}$ as
\begin{align}
\label{eqn:h_riemannian_metric}
     g^{(h,\tau)}(X,Y)\coloneqq\big\langle X,Y\big\rangle_{p}&\coloneq\sum_x X^{(m)}(\mathrm{x})Y^{(h,\tau)}(\mathrm{x})\nonumber \\
     &=\mathbb{E}_{p_{(h,\tau)}}\big[X^{(h,\tau)}Y^{(h,\tau)}\big].
\end{align}
For an arbitrary smooth real-valued function $f$ defined on $\mathcal{P}$, the differential $(df)_{p} \in T^*_{p}(\mathcal{P})$ maps a tangent vector $X \in T_{p}(\mathcal{P})$ to $X(f)$. The gradient vector $(\text{grad}\, f)_{p} \in T_{p}(\mathcal{P})$, associated with $(df)_{p}$, satisfies:
\begin{align}
\label{eqn:gradient-def}
\langle (\text{grad}\, f)_{p}, X \rangle_{p} = X(f) = (df)_{p}(X),\quad
(\text{grad}\, f)_{p} = \sum_{i,j} (g^{ij}(p))^{(h,\tau)}\, (\partial_jf)_p\, \partial_i, \nonumber
\end{align}where $(g^{ij}(p))^{(h,\tau)}$ denotes the $(i,j)^{\text{th}}$ entry of the inverse of the matrix $G_n^{(h,\tau)}(p)=\left[\langle \partial_i,\partial_j\rangle_p\right]$
and the squared norm of the differential is given by
\begin{equation}
\label{eqn:differential 2-form}
\| (df)_{p} \|_{p}^2 = \| (\text{grad}\, f)_{p} \|_{p}^2 = \sum_{i,j} (g^{ij}(p))^{(h,\tau)}\, (\partial_jf)_p\, (\partial_if)_p.
\end{equation}
The following theorem establishes a relationship between the variance of a random variable \(A\) under the escort distribution \(p_{(h,\tau)}\) and the geometric quantity \(\|(d\mathbb{E}[A])_p\|_p^2\).
\medskip
\begin{theorem}\label{thm:Variance_geo_theorem}
      For a random variable $A$ (a function from $\mathcal{X}$ to $\mathbb{R}$), let $\mathbb{E}[A]$ denote the mapping $p\mapsto \mathbb{E}_{p}[A]$ on $\mathcal{P}$. Then
\begin{align*}
    V_{p_{(h,\tau)}}[A]=\|(d\mathbb{E}[A])_{p}\|_{p}^{2},
\end{align*}
where the variance $ V_{p_{(h,\tau)}}[A]\coloneq \mathbb{E}_{p_{(h,\tau)}}[(A-\mathbb{E}_{p_{(h,\tau)}}[A])^2]$ and the norm $ \| \cdot \|_{p} $ is induced by the $(h,\tau)$-Fisher information metric.
\end{theorem}
\begin{proof}\
    For every $X\in T_{p}(\mathcal{P})$, we have
\[
X\big(\mathbb{E}[A]\big) = \big\langle(\text{grad}\,\mathbb{E}[A])_p, X\big\rangle_{p}.
\]
We observe that
\begin{align*}
     X\big(\mathbb{E}[A]\big)&=\sum_{x}X^{(m)}(x)A(x)\\
     &=\mathbb{E}_{p_{(h,\tau)}}[AX^{(h,\tau)}].
\end{align*} Since $\mathbb{E}_{p_{(h,\tau)}}[X^{(h,\tau)}]=0$, we can write
\begin{align}
     X\big(\mathbb{E}[A]\big)
     &=\mathbb{E}_{p_{(h,\tau)}}\left[X^{(h,\tau)}\left(A-\mathbb{E}_{p_{(h,\tau)}}[A]\right)\right].
\end{align} Observe that $\mathbb{E}_{p_{(h,\tau)}}\left\{A-\mathbb{E}_{p_{(h,\tau)}}[A]\right\}=0$. Thus, $A-\mathbb{E}_{p_{(h,\tau)}}[A]\in T_p^{(h,\tau)}(\mathcal{P})$, for some $Y\in T_p(\mathcal{P})$. Consequently, we can write $$X(\mathbb{E}[A])=\mathbb{E}_{p_{(h,\tau)}}[X^{(h,\tau)}Y^{(h,\tau)}].$$
Therefore,
\begin{align}
\label{eqn: metric expectation}  
\mathbb{E}_{p_{(h,\tau)}}[X^{(h,\tau)}Y^{(h,\tau)}]= \big\langle(\text{grad}\,\mathbb{E}[A])_p, X\big\rangle_{p}.
\end{align}
Again, from the definition of the $(h,\tau)$-Fisher information metric \eqref{eqn:h_riemannian_metric}, we see that 
\begin{equation}
\label{eqn:Ex-Innerproduct}
\mathbb{E}_{p_{(h,\tau)}}[X^{(h,\tau)}Y^{(h,\tau)}]= \big\langle Y, X\big\rangle_{p}, \quad \text{for any } X, Y \in T_{p}(\mathcal{P}).
\end{equation}
From \eqref{eqn: metric expectation} and \eqref{eqn:Ex-Innerproduct}, we get $Y = (\text{grad}\,\mathbb{E}[A])_p \in T_{p}(\mathcal{P}) $. Hence
\begin{align}
    \|(d\mathbb{E}[A])_{p}\|_{p}^{2}
    &= \big\langle(\operatorname{grad}\mathbb{E}[A])_p,
    (\operatorname{grad}\mathbb{E}[A])_p\big\rangle_{p}\nonumber\\
    &= \langle Y,Y\rangle_{p}\nonumber\\
    &= \mathbb{E}_{p_{(h,\tau)}}
    \left[X^{(h,\tau)}Y^{(h,\tau)}\right]\nonumber\\
    &= \mathbb{E}_{p_{(h,\tau)}}
    \left[\left(A-\mathbb{E}_{p_{(h,\tau)}}[A]\right)^2\right]\nonumber\\
    &=V_{p_{(h,\tau)}}[A].
    \label{eqn:star1}
\end{align}
\end{proof}
\begin{corollary}
  \label{corollarry 1}
 If $S$ is a sub-manifold of $\mathcal{P}$, then for $p\in S$ $$V_{p_{(h,\tau)}}[A]\geq\|(d\mathbb{E}[A]\mid_{S})_{p}\|_{p}^{2},$$ with equality if and only if
$$A-\mathbb{E}_{p_{(h,\tau)}}[A]\in T_{p}^{(h,\tau)}(S).$$  
\end{corollary}
\begin{proof}\ Since $(\operatorname{grad}f|_S)_p$ is the orthogonal projection of $(\operatorname{grad}f)_p$ onto $T_p(S)$, the result follows from Theorem~\ref{thm:Variance_geo_theorem}.
\end{proof}
We now define a differentiable function
$r:\Theta\to\mathbb{R}^k$ given by
\[
r(\theta)
=
\big[r_1(\theta),\ldots,r_k(\theta)\big]^\top,
\]
its Jacobian matrix is defined by
\[
J(\theta)
\coloneq\frac{\partial r(\theta)}{\partial\theta}\coloneqq
\big[\partial_j r_i(\theta)\big]_{i,j=1}^k.
\]
For an unbiased estimator
\[
\hat r(\mathrm{X})
=
\big[\hat r_1(\mathrm{X}),\ldots,\hat r_k(\mathrm{X})\big]^\top
\]
of $r(\theta)$, define its covariance matrix with respect to the $(h,\tau)$-escort distribution
$p_{(h,\tau,\theta)}$ by
\[V_{p_{(h,\tau,\theta)}}[\hat r(\mathrm{X})]\coloneqq
V_{(h,\tau,\theta)}[\hat r(\mathrm{X})]
\coloneq
\left[
\operatorname{Cov}_{(h,\tau,\theta)}
\big(\hat r_i(\mathrm{X}),\hat r_j(\mathrm{X})\big)
\right]_{i,j=1}^k,
\]
where
\[
\operatorname{Cov}_{(h,\tau,\theta)}
\big(\hat r_i(\mathrm{X}),\hat r_j(\mathrm{X})\big)
\coloneq
\mathbb{E}_{(h,\tau,\theta)}
\big[\hat r_i(\mathrm{X})\hat r_j({X})\big]
-
\mathbb{E}_{(h,\tau,\theta)}
\big[\hat r_i(\mathrm{X})\big]\,
\mathbb{E}_{(h,\tau,\theta)}
\big[\hat r_j(\mathrm{X})\big].
\]
Now, in the following, we state the main result of this section.
\medskip
\begin{theorem}[$(h,\tau)$-CR inequality]\label{thm:discrete_phi_crlb}
Let $\mathrm{X}=(X_1,\ldots,X_n)$ be distributed according to $p_\theta$. Let
\[
\hat r(\mathrm{X})
=
\big[\hat r_1(\mathrm{X}),\ldots,\hat r_k(\mathrm{X})\big]^\top
\]
be an unbiased estimator of $r(\theta)$ with respect to $p_\theta$, that is,
\[
\mathbb{E}_{\theta}[\hat r_i(\mathrm{X})]
=
r_i(\theta),
\qquad i=1,\ldots,k.
\]
Suppose that the $(h,\tau)$-Fisher information matrix
$G_n^{(h,\tau)}(\theta)$ is positive definite for every
$\theta\in\Theta$. Then
\begin{equation}
\label{eqn:discrete_h-CRLB}
V_{(h,\tau,\theta)}[\hat{r}(\mathrm{X})]
\succeq
J(\theta)
\left[G_n^{(h,\tau)}(\theta)\right]^{-1}
J(\theta)^\top,
\end{equation}where $G_n^{(h,\tau)}(\theta)=[g^{(h,\tau)}_{i,j}(\theta)]_{i,j=1}^k$ is the $(h,\tau)$-Fisher information matrix with $$g^{(h,\tau)}_{i,j}(\theta)=\mathbb{E}_{(h,\tau,\theta)}\left[\frac{1}{p_{(h,\tau,\theta)}(\mathrm{X})}\partial_ip_{\theta}(\mathrm{X})\frac{1}{p_{(h,\tau,\theta)}(\mathrm{X})}\partial_jp_{\theta}(\mathrm{X})\right].$$
\end{theorem}

\begin{proof}\
Define $A \coloneqq c^\top \hat{r}(\mathrm{X})$, where $\hat{r}(\mathrm{X})$ is an unbiased estimator of $r(\theta)$ with respect to $p_{\theta}$ and $c \in \mathbb{R}^k$. 
    Now for $p_{\theta}\in S$, we have
    \begin{align*}
         \|(d\mathbb{E}[A]\mid_{S})_{p_{\theta}}\|_{p_{\theta}}^2&\overset{(a)}{=}\frac{\partial}{\partial \theta}\mathbb{E}_{\theta}[c^{\top}\hat{r}(\mathrm{X})][G_{n}^{(h,\tau)}(\theta)]^{-1}\left\{\frac{\partial}{\partial \theta}\mathbb{E}_{\theta}[c^{\top}\hat{r}(\mathrm{X})]\right\}^{\top},\\
         &=c^{\top} \frac{\partial \mathbb{E}_{\theta}[\hat{r}(\mathrm{X})]}{\partial \theta}[G_n^{(h,\tau)}(\theta)]^{-1}\left\{\frac{\partial \mathbb{E}_{\theta}[\hat{r}(\mathrm{X})]}{\partial \theta}\right\}^{\top}c\\
         &=c^{\top}J(\theta)[G_n^{(h,\tau)}(\theta)]^{-1} J(\theta)^{\top} c,
    \end{align*}where $(a)$ follows from \eqref{eqn:differential 2-form},
    and
    \begin{equation*}
   V_{(h,\tau,\theta)}[A]=c^{\top}V_{(h,\tau,\theta)}[\hat{r}(\mathrm{X})]c.
    \end{equation*}\\ 
    Now applying the Corollary \ref{corollarry 1} to the submanifold $S$, we get
 $$c^{\top} V_{(h,\tau,\theta)}[\hat{r}(\mathrm{X})]\;c\geq c^{\top}J(\theta)[G_n^{(h,\tau)}(\theta)]^{-1}J(\theta)^\top c.$$
This establishes \eqref{eqn:discrete_h-CRLB}.
\end{proof}
In the following, we present the necessary and sufficient condition for $(h,\tau)$-CR inequality.
\begin{proposition}\label{prop:discrete equality condition}
    Equality holds in \eqref{eqn:discrete_h-CRLB} for $\hat{r}(\mathrm{X})$ if and only if
\begin{equation}\label{eqn:equality phi CRLB_finite support}
  \hat{r}(\mathrm{x})-\mathbb{E}_{(h,\tau,\theta)}[\hat{r}(\mathrm{X})]=J(\theta)\big[G_n^{(h,\tau)}(\theta)\big]^{-1}\frac{1}{p_{(h,\tau,\theta)}(\mathrm{x})}\frac{\partial p_{\theta}(\mathrm{x})}{\partial \theta}.
\end{equation} 
\end{proposition}
\begin{proof}\
From Corollary~\ref{corollarry 1}, we get that equality in \eqref{eqn:discrete_h-CRLB} if and only if $$\sum_{i=1}^k c_i(\hat{r}_i(\mathrm{x})-\mathbb{E}_{(h,\tau,\theta)}[\hat{r}_i(\mathrm{X})])\in T_{p_{\theta}}^{(h,\tau)}(S)\text{ for all }[c_i]\in\mathbb{R}^k.$$
That is,
\begin{align*}
   &\mathbb{E}_{(h,\tau,\theta)}\left[\sum_{i=1}^k c_i\left(\hat{r}_i(\mathrm{X})-\mathbb{E}_{(h,\tau,\theta)}[\hat{r}_i(\mathrm{X})]\right)\right]=0\text{ for all } [c_i]\in \mathbb{R}^k.
\end{align*}
This implies that 
$$c^\top \mathbb{E}_{(h,\tau,\theta)}\bigg[\hat{r}(\mathrm{X})-\mathbb{E}_{(h,\tau,\theta)}[\hat{r}(\mathrm{X})]\bigg]=0.$$
The above is true for every $c\in \mathbb{R}^k$.\\ 
Hence
\begin{align*}
  \mathbb{E}_{(h,\tau,\theta)}\bigg[\hat{r}_i(\mathrm{X})-\mathbb{E}_{(h,\tau,\theta)}[\hat{r}_i(\mathrm{X})]\bigg]=0\text{ for all }i.
\end{align*}That is,
$$\hat{r}_i(\mathrm{x})-\mathbb{E}_{(h,\tau,\theta)}[\hat{r}_i(\mathrm{X})]\in T_{p_{\theta}}^{(h,\tau)}(S) \text{ for all }i.$$
This implies that
\begin{align*}
 \hat{r}(\mathrm{x})-\mathbb{E}_{(h,\tau,\theta)}[\hat{r}(\mathrm{X})]
   =&
   \begin{bmatrix}
       \hat{r}_1(\mathrm{x})-\mathbb{E}_{(h,\tau,\theta)}[\hat{r}_1(\mathrm{X})]\\ 
       \vdots\\
      \hat{r}_k(\mathrm{x})-\mathbb{E}_{(h,\tau,\theta)}[\hat{r}_k(\mathrm{X})]
   \end{bmatrix}\\
   =&\begin{bmatrix}
       \sum_{j=1}^k m_{1j}(\theta)\frac{1}{p_{(h,\tau,\theta)}(\mathrm{x})}\frac{\partial p_{\theta}(\mathrm{x})}{\partial \theta_j} \\
       \vdots\\
       \sum_{j=1}^k m_{kj}(\theta)\frac{1}{p_{(h,\tau,\theta)}(\mathrm{x})}\frac{\partial p_{\theta}(\mathrm{x})}{\partial \theta_j}
   \end{bmatrix}\\
    =&\begin{bmatrix}
        m_{11}(\theta)&\hdots &m_{1k}(\theta)\\
        \vdots &\vdots &\vdots\\
        m_{k1}(\theta)& \hdots &m_{kk}(\theta)\\
    \end{bmatrix}
    \begin{bmatrix}
 \frac{1}{p_{(h,\tau,\theta)}(\mathrm{x})}\frac{\partial p_{\theta}(\mathrm{x})}{\partial \theta_1} \\
        \vdots\\
     \frac{1}{p_{(h,\tau,\theta)}(\mathrm{x})}\frac{\partial p_{\theta}(\mathrm{x})}{\partial \theta_k} 
    \end{bmatrix}\\
    &\;\text{ for some } m_{ij}(\theta), \text{ functions of }\theta,\; \text{where } i,j=1,\dots,k\\
    =&\, m(\theta)\frac{1}{p_{(h,\tau,\theta)}({x})}\frac{\partial p_{\theta}(\mathrm{x})}{\partial \theta}.
\end{align*} That is, 
\begin{equation}
\label{eqn:star2}
     \hat{r}(\mathrm{x})-\mathbb{E}_{(h,\tau,\theta)}[\hat{r}(\mathrm{X})]=m(\theta)\frac{1}{p_{(h,\tau,\theta)}(\mathrm{x})}\frac{\partial p_{\theta}(\mathrm{x})}{\partial \theta} ,
\end{equation}
for some $k\times k$ matrix $m(\theta)$.\\
We now determine $m(\theta)$. Multiplying $u_{\theta}(\mathrm{x})^{\top}$ on the both sides of \eqref{eqn:star2}, we get 
\begin{equation*}
    \hat{r}(\mathrm{x})u_{\theta}(\mathrm{x})^{\top}-\mathbb{E}_{(h,\tau,\theta)}[\hat{r}(\mathrm{X})]u_{\theta}(\mathrm{x})^{\top}=m(\theta)\frac{1}{p_{(h,\tau,\theta)}({x})}\frac{\partial p_{\theta}(\mathrm{x})}{\partial\theta}u_{\theta}(\mathrm{x})^{\top}.
\end{equation*}Taking expectation with respect to $p_{\theta}$, we obtain
\begin{align*}
      \sum_{\mathrm{x}}\hat{r}(\mathrm{x})u_{\theta}(\mathrm{x})^{\top}p_{\theta}(\mathrm{x})-\mathbb{E}_{(h,\tau,\theta)}[\hat{r}(X)]&\sum_{\mathrm{x}} u_{\theta}(\mathrm{x})^{\top}p_{\theta}({x})\\
      &=m(\theta)\mathbb{E}_{(h,\tau,\theta)}\left[\frac{1}{p_{(h,\tau,\theta)}({X})}\frac{\partial}{\partial\theta}p_{\theta}({X})\frac{1}{p_{(h,\tau,\theta)}({X})}\frac{\partial}{\partial\theta}p_{\theta}({X})^{\top}\right].
\end{align*}
Since $\mathbb{E}_{\theta}[u_{\theta}(\mathrm{X})]=0$, the second term on the left-hand side vanishes. Therefore, 
\begin{equation*}
    \frac{\partial \mathbb{E}_{\theta}[\hat{r}(\mathrm{X})]}{\partial \theta}=m(\theta)G_{n}^{(h,\tau)}(\theta).
\end{equation*}Hence
\begin{equation*}
   J(\theta)=m(\theta)G_{n}^{(h,\tau)}(\theta),
\end{equation*}
which implies, \begin{equation*}
    m(\theta)= J(\theta)[G_{n}^{(h,\tau)}(\theta)]^{-1}.
\end{equation*}
Substituting this in \eqref{eqn:star2}, we obtain the equality condition \eqref{eqn:equality phi CRLB_finite support}.
\end{proof} 

\begin{remark}
    For $h(t)=t\log t-t+1$, $\tau(t)=t$, \eqref{eqn:discrete_h-CRLB} coincides with the classical CR inequality \eqref{eqn:CRLB}. Furthermore, the equality condition \eqref{eqn:equality phi CRLB_finite support} coincides with  \begin{equation}\label{eqn:CRLB_equality_condition}
       \hat{r}(\mathrm{x})-\mathbb{E}_{\theta}[\hat{r}(\mathrm{X})]=J(\theta)[G_n(\theta)]^{-1}u_{\theta}(\mathrm{x}),
   \end{equation}which is the equality condition for \eqref{eqn:CRLB}.
\end{remark}
\begin{remark}
For $h(t)=\varphi(t)$, $\tau(t)=t$, inequality
\eqref{eqn:discrete_h-CRLB} coincides with the
$\varphi$-CRLB established in \cite{dhadumia2026cramer}.
\end{remark}
We present the following example to illustrate \eqref{eqn:discrete_h-CRLB}. 

 \begin{example}\label{example:inequality_equality}
     Let $\mathrm{X}=(X_1,X_2)$ be distributed according to 
     \begin{equation}\label{exp:Modi_Ber}
         p_{\theta}(\mathrm{x})=\frac{(1-\theta)}{2}\left[1+\frac{(2\theta-1)}{(1-\theta)}\frac{(x_1+x_2)}{2}\right], \quad\mathrm{x}=(x_1,x_2)\in \{0,1\}^2, \quad\theta\in\left(1/2,1\right).
     \end{equation}
     We can express \eqref{exp:Modi_Ber} as follows
     \begin{equation}\label{eqn:Ber_alpha_form}
         p_{\theta}(\mathrm{x})=\frac{3-2\theta}{8}+\frac{2\theta-1}{4}\left(x_1+x_2-\frac{1}{2}\right).
     \end{equation}
     Consider the unbiased estimator $T(\mathrm{X})=X_1+X_2-\frac{1}{2}$ of $\theta$. We consider the following two choices of $h$ and $\tau$ and evaluate the corresponding $(h,\tau)$-CR inequalities.
     \begin{itemize}
  
\item[(1)] Choose $h(t)=\frac{2}{3}t^{3/2}$ and
$\tau(t)=t^2$. Then
\[
p_{(h,\tau,\theta)}(\mathrm{x})
=
\frac{p_\theta(\mathrm{x})^{-1}}
{\sum_{\mathrm{y}}p_\theta(\mathrm{y})^{-1}},
\]
where
\[
\begin{aligned}
p_{(h,\tau,\theta)}(0,0)
&=\frac{\theta}{1+4\theta(1-\theta)},\\
p_{(h,\tau,\theta)}(0,1)
&=p_{(h,\tau,\theta)}(1,0)
=\frac{2\theta(1-\theta)}{1+4\theta(1-\theta)},\\
p_{(h,\tau,\theta)}(1,1)
&=\frac{(1-\theta)}{1+4\theta(1-\theta)}.
\end{aligned}
\]
We compute
\[
\begin{aligned}
V_{(h,\tau,\theta)}[T(\mathrm{X})]
&=\frac{8\theta(1-\theta)}
{[1+4\theta(1-\theta)]^2},\\
J(\theta)[G_n^{(h,\tau)}(\theta)]^{-1}J(\theta)
&=\frac{4\theta(1-\theta)}
{1+4\theta(1-\theta)}.
\end{aligned}
\]
For $\theta\in\left(\frac{1}{2},1\right)$,
\begin{equation}
V_{(h,\tau,\theta)}[T(\mathrm{X})]
>
J(\theta)[G_n^{(h,\tau)}(\theta)]^{-1}J(\theta).\label{eqn:st_inequality}    
\end{equation}
\item[(2)] Choose $h(t)=\frac{1}{2}t^2$ and $\tau(t)=t$.
Then $p_{(h,\tau,\theta)}(\mathrm{x})=\frac{1}{4}$ for all $\mathrm{x}\in \{0,1\}^2$.
In this case, equality holds:
\[
V_{(h,\tau,\theta)}[T(\mathrm{X})]
=
J(\theta)[G_n^{(h,\tau)}(\theta)]^{-1}J(\theta)
=
\frac{1}{2}.
\]Observe that $\hat{r}(\mathrm{x})-\mathbb{E}_{(h,\tau,\theta)}[\hat{r}(\mathrm{X})]=J(\theta)\big[G_n^{(h,\tau)}(\theta)\big]^{-1}\frac{1}{p_{(h,\tau,\theta)}(\mathrm{x})}\frac{\partial p_{\theta}(\mathrm{x})}{\partial \theta}=x_1+x_2-1$.
     \end{itemize}
 \end{example}
\section{{Extension to General Support Sets}}\label{sec:3}
  Due to the limitation of Amari and Nagaoka's theory to models with finite support, \eqref{eqn:discrete_h-CRLB} is applicable only to such models. In the following, we extend this inequality to models with more general support sets under suitable regularity assumptions, using Cauchy–Schwarz inequality.
  
  Let $\mathrm{X}=(X_1, \dots, X_n)$ be distributed according to $p_{\theta}$ of support $\mathcal{X}_{\theta}$. Let $\hat{r}(\mathrm{X})=[\hat{r}_1(\mathrm{X}),\ldots,\hat{r}_k(\mathrm{X})]^{\top}$ be an unbiased estimator of the parameter $r(\theta)=[r_1(\theta),\ldots,r_k(\theta)]^{\top}$, a $k$-dimensional differentiable function of $\theta$, with respect to $p_{\theta}$, that is, $\mathbb{E}_{\theta}[\hat{r}_i(\mathrm{X})]=r_i(\theta)$ for $i=1,\ldots,k$. We assume the following regularity conditions:
\begin{itemize}
    \item [(I)] $\int p_{\theta}(\mathrm{x})d\mu(\mathrm{x})$ can be differentiated under the integral sign with respect to $\theta$, that is,
  $$
 \frac{\partial}{\partial\theta} \int p_{\theta}(\mathrm{x}) \, d\mu(\mathrm{x})= \int \frac{\partial}{\partial\theta} p_{\theta}(\mathrm{x}) d\mu(\mathrm{x}).
  $$ \\
  \item [(II)]$\int \hat{r}(\mathrm{x})p_{\theta}(\mathrm{x})d\mu(\mathrm{x})$ is finite and can be differentiated under the integral sign with respect to $\theta$, that is, 
   $$\frac{\partial}{\partial \theta}\int \hat{r}(\mathrm{x})p_{\boldsymbol\theta}(\mathrm{x})d\mu(\mathrm{x})=\int \hat{r}(\mathrm{x})\frac{\partial}{\partial\theta} p_{\theta}(\mathrm{x})^{\top} d\mu(\mathrm{x}).$$
   \item[(III)] The functions $ \mathbb{X}(\mathrm{x}) \coloneqq \hat{r}(\mathrm{x}) - \mathbb{E}_{(h,\tau, \theta)}[\hat{r}(\mathrm{X})]$ and $\mathbb{Y}(\mathrm{x}) \coloneqq \frac{1}{p_{(h,\tau,\theta)}(\mathrm{x})}\frac{\partial}{\partial \theta}p_{\theta}(\mathrm{x})$ are square-integrable with respect to the $(h,\tau)$-escort distribution $$p_{(h,\tau,\theta)}(\mathrm{x})\coloneqq\frac{[h''(\tau(p_{\theta}(\mathrm{x})))[\tau'(p_{\theta}(\mathrm{x}))]^2]^{-1}}{\int[h''(\tau(p_{\theta}(\mathrm{y})))[\tau'(p_{\theta}(\mathrm{y}))]^2]^{-1}d\mu(\mathrm{y})}.$$\\
   \item[(IV)] The $(h,\tau)$-Fisher information matrix $G^{(h,\tau)}_n(\theta)\coloneqq \mathbb{E}_{(h,\tau,\theta)}\left[\frac{1}{p_{(h,\tau,\theta)}(\mathrm{X})}\frac{\partial p_{\theta}(\mathrm{X})}{\partial\theta}\frac{1}{p_{(h,\tau,\theta)}(\mathrm{X})}\frac{\partial p_{\theta}(\mathrm{X})}{\partial\theta}^{\top}\right]$ is positive definite for every $\theta\in\Theta$.
\end{itemize}
\begin{theorem} \label{thm:h_CRLB} Under the assumptions in conditions $\text{(I)}-\text{(IV)}$, we have \begin{equation}
\label{eqn:h-CRLB}
V_{(h,\tau,\theta)}[\hat{r}(\mathrm{X})]
\succeq
J(\theta)
\left[G_n^{(h,\tau)}(\theta)\right]^{-1}
J(\theta)^\top.
\end{equation}Furthermore, equality holds in \eqref{eqn:h-CRLB} if and only if 
\begin{equation}\label{eqn:equality_condition_h_crlb}
     \hat{r}(\mathrm{x})-\mathbb{E}_{(h,\tau,\theta)}[\hat{r}(\mathrm{X})]=J(\theta)\big[G_n^{(h,\tau)}(\theta)\big]^{-1}\frac{1}{p_{(h,\tau,\theta)}(\mathrm{x})}\frac{\partial p_{\theta}(\mathrm{x})}{\partial \theta}.
\end{equation}
    
\end{theorem}
\begin{proof}\
    For any $a \in \mathbb{R}^k$, let $\mathbf{f}(\mathrm{x})=a^{\top}\mathbb{X}(\mathrm{x})$ and $\mathbf{g}(\mathrm{x})=\mathbb{Y}(\mathrm{x})^{\top}b$, where $b=[G_n^{(h,\tau)}(\theta)]^{-1}
J(\theta)^{\top} a$. Then $\mathbf{f}$ and $\mathbf{g}$ are square-integrable with respect to $p_{(h,\tau,\theta)}$ by $\text{III}$. Now, applying Cauchy–Schwarz inequality to the functions $\mathbf{f}$ and $\mathbf{g}$, we get
\begin{equation}
    \label{Cauchy-schwarz ine}
    \big(\mathbb{E}_{(h,\tau,\theta)}\left[\mathbf{f}(\mathrm{X})\cdot \mathbf{g}(\mathrm{X}) \right]\big)^2
\leq 
\mathbb{E}_{(h,\tau,\theta)}\left[\mathbf{f}(\mathrm{X})^2\right] \cdot 
\mathbb{E}_{(h,\tau,\theta)}\left[\mathbf{g}(\mathrm{X}) ^2  \right].
\end{equation}
Observe that
\begin{align*}
   \mathbb{E}_{(h,\tau,\theta)}&[\mathbf{f}(\mathrm{X})\cdot \mathbf{g} (\mathrm{X})]=\int a^{\top}\big(\hat{r}(\mathrm{x})-\mathbb{E}_{(h,\tau,\theta)}[\hat{r}(\mathrm{X})]\big)\frac{1}{p_{(h,\tau,\theta)}(\mathrm{x})}\frac{\partial}{\partial \theta}p_{\theta}(\mathrm{x})^{\top}b\,p_{(h,\tau,\theta)}(\mathrm{x})d\mathrm{x}\\
   =& a^{\top}\int\big(\hat{r}(\mathrm{x})-\mathbb{E}_{(h,\tau,\theta)}[\hat{r}(\mathrm{X})]\big)\frac{\partial}{\partial \theta}p_{\theta}(\mathrm{x})^{\top}d\mathrm{x}\cdot b\\
    =&a^{\top}\left[\int\hat{r}(\mathrm{x})\frac{\partial}{\partial\theta}p_{\theta}(\mathrm{x})^{\top}d\mathrm{x}-\int \mathbb{E}_{(h,\tau,\theta)}[\hat{r}(\mathrm{X})]\frac{\partial}{\partial\theta}p_{\theta}(\mathrm{x})^{\top}d\mathrm{x}\right]b\\
    \overset{(a)}{=}& a^{\top}\left\{ \frac{\partial \mathbb{E}_{\theta}[\hat{r}(\mathrm{X})]}{\partial \theta}- \mathbb{E}_{(h,\tau,\theta)}[\hat{r}(\mathrm{X})]\left(\frac{\partial}{\partial\theta}\int  p_{\theta}(\mathrm{x})d\mathrm{x}\right)^{\top}\right\}b\\
    =&a^{\top}J(\theta)\;b,
   \end{align*}
where equality $(a)$ follows from regularity conditions $\text{I}$ and $\text{II}$.
\\ \\
$\mathbb{E}_{(h,\tau,\theta)}\left[\mathbf{f}(\mathrm{X})^2 \right]=a^{\top}\mathbb{E}_{(h,\tau,\theta)}[(\hat{r}(\mathrm{X})-\mathbb{E}_{(h,\tau,\theta)}[\hat{r}(\mathrm{X})])(\hat{r}(\mathrm{X})-\mathbb{E}_{(h,\tau,\theta)}[\hat{r}(\mathrm{X})])^{\top}]a
    =a^{\top}V_{(h,\tau,\theta)}[\hat{r}(\mathrm{X})]a,$ which is finite since $\mathbf{f}$ is square-integrable.
\begin{align*}
    \mathbb{E}_{(h,\tau,\theta)}[\mathbf{g}(\mathrm{X})^2]&=b^{\top}\mathbb{E}_{(h,\tau,\theta)}\left[\frac{1}{p_{(h,\tau,\theta)}(\mathrm{X})}\frac{\partial}{\partial\theta}p_{\theta}(\mathrm{X})\frac{1}{p_{(h,\tau,\theta)}(\mathrm{X})}\frac{\partial}{\partial\theta}p_{\theta}(\mathrm{X})^{\top}\right]b\\
   &=b^{\top}G_{n}^{(h,\tau)}(\theta)b.
\end{align*}
From \eqref{Cauchy-schwarz ine}, we get
\begin{align}
    \bigg[a^\top  J(\theta)b\bigg]^2\leq a^{\top} V_{(h,\tau,\theta)}[\hat{r}(\mathrm{X})]a b^{\top}G_{n}^{(h,\tau)}(\theta)b.\label{Cauchy-Swartz in variance terms}
\end{align}
Since $b=[G_n^{(h,\tau)}(\theta)]^{-1}J(\theta)^{\top}a$, we get
\begin{align}
\label{Cauchy-Swartz in variance terms LHS}
   \bigg[a^\top J(\theta)b\bigg]^2 =\bigg[a^\top J(\theta)[G_n^{(h,\tau)}(\theta)]^{-1}J(\theta)^{\top}a\bigg]^2
\end{align}
and
\begin{align}
\label{Cauchy-Swartz in variance terms RHS}
    a^{\top} V_{(h,\tau,\theta)}[\hat{r}(\mathrm{X})]ab^{\top}G_{n}^{(h,\tau)}(\theta)b
    &=a^{\top} V_{(h,\tau,\theta)}[\hat{r}(\mathrm{X})]a a^{\top}J(\theta)[G_n^{(h,\tau)}(\theta)]^{-1}G_n^{(h,\tau)}(\theta)[G_n^{(h,\tau)}(\theta)]^{-1}J(\theta)^{\top}a\nonumber\\
    &=a^{\top} V_{(h,\tau,\theta)}[\hat{r}(\mathrm{X})]a \Bigg[a^{\top}J(\theta)[G_n^{(h,\tau)}(\theta)]^{-1}
   J(\theta)^{\top}a\Bigg].
\end{align} Since $G_n^{(h,\tau)}(\theta)$ is positive definite, so is $[G_n^{(h,\tau)}(\theta)]^{-1}$. Therefore $J(\theta)[G_n^{(h,\tau)}(\theta)]^{-1}
   J(\theta)^{\top}$ is at least positive semi-definite.
Hence, using \eqref{Cauchy-Swartz in variance terms LHS} and \eqref{Cauchy-Swartz in variance terms RHS} in \eqref{Cauchy-Swartz in variance terms}, we get 
\begin{equation}
    a^{\top}J(\theta)[G_n^{(h,\tau)}(\theta)]^{-1}
   J(\theta)^{\top}a\leq a^{\top}V_{(h,\tau,\theta)}[\hat{r}(\mathrm{X})]a.
\end{equation}
This establishes $(h,\tau)$-CR inequality. Furthermore, by the equality condition of Cauchy–Schwarz inequality \eqref{Cauchy-schwarz ine}, we get equality in \eqref{eqn:h-CRLB} if and only if
\begin{align}
a^{\top}\left\{\hat{r}(\mathrm{x}) - \mathbb{E}_{(h,\tau,\theta)}[\hat{r}(\mathrm{X})]\right\}
= c(\theta)\, \frac{1}{p_{(h,\tau,\theta)}(\mathrm{x})}\frac{\partial}{\partial\theta}p_{\theta}(\mathrm{x})^{\top}b, \nonumber
\end{align}
where $c(\theta)$ is a constant independent of $\mathrm{x}$.
That is, 
\begin{align}
a^{\top}\left\{\hat{r}(\mathrm{x}) - \mathbb{E}_{(h,\tau,\theta)}[\hat{r}(\mathrm{X})]\right\}
&=  c(\theta)\, \frac{1}{p_{(h,\tau,\theta)}(\mathrm{x})}\frac{\partial}{\partial\theta}p_{\theta}(\mathrm{x})^{\top}[G_n^{(h,\tau)}(\theta)]^{-1}J(\theta)^{\top}a\nonumber\\
&=a^{\top}c(\theta)J(\theta)[G_n^{(h,\tau)}(\theta)]^{-1}\frac{1}{p_{(h,\tau,\theta)}(\mathrm{x})}\frac{\partial}{\partial\theta}p_{\theta}(\mathrm{x}),\label{eqn:eq_cond_intermed}
\end{align}
where \eqref{eqn:eq_cond_intermed} is obtained by using the fact that  $s^\top = s$ for scalar $s$. Since $a$ is arbitrary, \eqref{eqn:eq_cond_intermed} gives 
\begin{equation}\label{cond:Cauchy_swartz_eq}
 \hat{r}(\mathrm{x}) - \mathbb{E}_{(h,\tau,\theta)}[\hat{r}(\mathrm{X})]  =c(\theta)J(\theta)[G_n^{(h,\tau)}(\theta)]^{-1}\frac{1}{p_{(h,\tau,\theta)}(\mathrm{x})}\frac{\partial}{\partial\theta}p_{\theta}(\mathrm{x}).
\end{equation}
Multiplying \eqref{cond:Cauchy_swartz_eq} by $u_{\theta}(\mathrm{x})^{\top}$, we get
\begin{equation*}
    \hat{r}(\mathrm{x})u_{\theta}(\mathrm{x})^{\top} - \mathbb{E}_{(h,\tau,\theta)}[\hat{r}(\mathrm{X})]  u_{\theta}(\mathrm{x})^{\top}=c(\theta)J(\theta)[G_n^{(h,\tau)}(\theta)]^{-1}\frac{1}{p_{(h,\tau,\theta)}(\mathrm{x})}\frac{\partial}{\partial\theta}p_{\theta}(\mathrm{x})u_{\theta}(\mathrm{x})^{\top}.
\end{equation*}
Taking expectation with respect to $p_{\theta}$, we obtain
\begin{align*}
    \int   \hat{r}(\mathrm{x})u_{\theta}(\mathrm{x})^{\top}p_{\theta}(\mathrm{x})d\mathrm{x}-&\mathbb{E}_{(h,\tau,\theta)}[\hat{r}(\mathrm{X})]  \int u_{\theta}(\mathrm{x})^{\top}p_{\theta}(\mathrm{x})d\mathrm{x}\\
    =&c(\theta)J(\theta)[G_n^{(h,\tau)}(\theta)]^{-1}\int \frac{1}{p_{(h,\tau,\theta)}(\mathrm{x})}\frac{\partial}{\partial\theta}p_{\theta}(\mathrm{x})\frac{\partial}{\partial\theta}p_{\theta}(\mathrm{x})^{\top}\\
    =&c(\theta)J(\theta)[G_n^{(h,\tau)}(\theta)]^{-1}G_n^{(h,\tau)}(\theta).
\end{align*}This implies
\begin{equation*}
    \frac{\partial \mathbb{E}_{\theta}[\hat{r}(\mathrm{X})]}{\partial \theta}=c(\theta)J(\theta).
\end{equation*}
Hence \begin{equation*}
 J(\theta)=c(\theta)J(\theta).
\end{equation*}
Substituting this in \eqref{cond:Cauchy_swartz_eq}, we obtain
\begin{equation*}
   \hat{r}(\mathrm{x})-\mathbb{E}_{(h,\tau,\theta)}[\hat{r}(\mathrm{X})]=J(\theta)\big[G_n^{(h,\tau)}(\theta)\big]^{-1}\frac{1}{p_{(h,\tau,\theta)}(\mathrm{x})}\frac{\partial}{\partial \theta}p_{\theta}(\mathrm{x}),
\end{equation*} which is the equality condition \eqref{eqn:equality_condition_h_crlb}. This completes the proof.
\end{proof} 

\section{A Generalized CR Inequality for $(h,\tau)$-Exponential Family}\label{sec:4}
In this section, we define the $(h,\tau)$-family of probability distributions, which is a generalization of the exponential family \cite[Eq.~4.1]{MATSUZOE2027131088}. It is well known that, under suitable regularity conditions, the equality in the classical CR inequality holds if and only if the underlying family of distributions is an exponential family \cite{wijsman1973attainment,fend1959attainment,muller1989attainment}. However, this attainment property does not, in general, extend to the $(h,\tau)$-CR inequality for an $(h,\tau)$-exponential family. Consequently, we establish a CR inequality whose lower bound is attained under an $(h,\tau)$-exponential family.
\begin{definition}[$(h,\tau)$-exponential family]
    Let $h$ and $\tau$ as in \eqref{eqn:(h,tau)-divergence}. A family of probability distributions $S = \{p_{\theta}:\theta\in\Theta\}$, where $\Theta$ is an open subset of $\mathbb{R}^k$, is said to form an $(h,\tau)$-exponential family if each $p_{\theta}\in S$ can be written as 
\begin{equation}
\label{eqn:h_tau_family}
(h'\circ\tau)\bigl(p_\theta(\mathrm{x})\bigr) =
 C(\mathrm{x}) + Z(\theta)+\lambda(\theta)^{\top}T(\mathrm{x})\qquad
\text{for }\mathrm{x}\in\mathcal{X},
\end{equation}
for some $C:\mathcal{X}\to\mathbb{R}$ and $Z\in C^1(\Theta)$, and $\lambda=[\lambda_1,\ldots,\lambda_k]^{\top}$, $T=[T_1,\ldots,T_k]^{\top}$, where $\lambda_i\in C^1(\Theta)$ and $T_i:\mathcal{X}\to \mathbb{R}$, $i=1,\ldots,k$.
\end{definition}
\begin{example}
    With $h(t)=\varphi(t)$, $\tau(t)=t$, $(h,\tau)$-exponential family reduces to $B^{(\varphi)}$-family in \cite{dhadumia2026cramer}.
\end{example}
\begin{example}
    With $h(t)=t\log t-t+1$, $\tau(t)=t$, $(h,\tau)$-exponential family reduces to exponential family.
\end{example}
\begin{remark}Theorem~6 of \cite{dhadumia2025generalized} establishes that equality in the $\varphi$-CRLB is attained under the $B^{(\varphi)}$-family. For the particular choice $(h(t),\tau(t))=(\varphi(t),t)$, the $(h,\tau)$-CR inequality and the $(h,\tau)$-exponential family reduce to the $\varphi$-CRLB and the $B^{(\varphi)}$-family, respectively. Hence, equality in the $(h,\tau)$-CR inequality is attained under the $(h,\tau)$-exponential family for this particular choice of $(h,\tau)$.
\end{remark}
\begin{remark}\label{remark:alpha_family_ofBer}
For $C(\mathrm{x})=0$, $Z(\theta)=\frac{3-2\theta}{8}$,
$\lambda(\theta)=\frac{2\theta-1}{4}$, and
$T(\mathrm{x})=x_1+x_2-\frac12$, the $(h,\tau)$-exponential
family in \eqref{eqn:h_tau_family} coincides with
\eqref{eqn:Ber_alpha_form} for both choices
$(h(t),\tau(t))=(\frac23t^{3/2},t^2)$ and
$(h(t),\tau(t))=(\frac12t^2,t)$. Hence, \eqref{eqn:Ber_alpha_form} forms an $(h,\tau)$-exponential family. 
However, as shown in Example~\ref{example:inequality_equality},
the $(h,\tau)$-CR inequality is strict for $T(\mathrm{X})$
in the first case, whereas equality holds in the second. Thus, in general, equality in the $(h,\tau)$-CR inequality
is not necessarily attained under the $(h,\tau)$-exponential family.  
\end{remark} 
In the following, we establish a CR inequality for which equality is attained under the \((h,\tau)\)-exponential family.\\

Let ${X}=(X_1, \dots, X_n)$ be distributed according to $p_{\theta}$ of support $\mathcal{X}_{\theta}$. Let $\hat r(\mathrm X)=[\hat r_1(\mathrm X),\ldots,\hat r_k(\mathrm X)]^\top$ be an unbiased estimator of $r(\theta)\coloneq[r_1(\theta),\ldots,r_k(\theta)]^{\top}$, a $k$-dimensional differentiable function. We now define $\widetilde{(h,\tau)}$-escort distribution $$\widetilde{p}_{(h,\tau,\theta)}(\mathrm{x})\coloneqq\frac{[h''(\tau(p_{\theta}(\mathrm{x})))\tau'(p_{\theta}(\mathrm{x}))]^{-1}}{\int [h''(\tau(p_{\theta}(\mathrm{y})))\tau'(p_{\theta}(\mathrm{y}))]^{-1}d\mathrm{y}}.$$ 
For the estimator
$\hat r(\mathrm X)=[\hat r_1(\mathrm X),\ldots,\hat r_k(\mathrm X)]^\top$, its covariance matrix under $\widetilde{p}_{(h,\tau,\theta)}$ is defined as 
\[
\widetilde{V}_{(h,\tau,\theta)}[\hat r(\mathrm{X})]
\coloneq
\left[
\widetilde{\operatorname{Cov}}_{(h,\tau,\theta)}
\big(\hat r_i(\mathrm{X}),\hat r_j(\mathrm{X})\big)
\right]_{i,j=1}^k,
\]
where
\[
\widetilde{\operatorname{Cov}}_{(h,\tau,\theta)}
\big(\hat r_i(\mathrm{X}),\hat r_j(\mathrm{X})\big)
\coloneq
\widetilde{\mathbb{E}}_{(h,\tau,\theta)}
\big[\hat r_i(\mathrm{X})\hat r_j({X})\big]
-
\widetilde{\mathbb{E}}_{(h,\tau,\theta)}
\big[\hat r_i(\mathrm{X})\big]\,
\widetilde{\mathbb{E}}_{(h,\tau,\theta)}
\big[\hat r_j(\mathrm{X})\big],
\] where $\widetilde{\mathbb{E}}_{(h,\tau,\theta)}$ is with respect to $\widetilde{p}_{(h,\tau,\theta)}$.

We now consider the following regularity conditions:
\begin{itemize}
    \item [(A)] The functions $X(\mathrm{x})\coloneqq \hat{r}(\mathrm{x})-\widetilde{\mathbb{E}}_{(h,\tau,\theta)}[\hat{r}(\mathrm{X})]$, and $Y(\mathrm{x})\coloneqq \frac{1}{\widetilde{p}_{(h,\tau,\theta)}(\mathrm{x})}\frac{\partial p_{\theta}(\mathrm{x})}{\partial\theta}$ are square integrable with respect to $\widetilde{p}_{(h,\tau,\theta)}$. \vspace{.2cm}
    \item[(B)] The $\widetilde{(h,\tau)}$-Fisher information matrix $$\widetilde{G}_n^{(h,\tau)}(\theta)\coloneq\widetilde{\mathbb{E}}_{(h,\tau,\theta)}\left[\frac{1}{\widetilde{p}_{(h,\tau,\theta)}(\mathrm{X})}\frac{\partial p_{\theta}(\mathrm{X})}{\partial\theta}\frac{1}{\widetilde{p}_{(h,\tau,\theta)}(\mathrm{X})}\frac{\partial p_{\theta}(\mathrm{X})}{\partial\theta}^{\top}\right]$$ 
\end{itemize} is positive definite for every $\theta\in\Theta$.
\begin{theorem}[$\widetilde{(h,\tau)}$-CR inequality]\label{thm:t_CRLB} Under the assumptions in conditions $(I)-(II)$ and $(A)-(B)$, we have
   \begin{equation}\label{eqn:t_CRLB}
       \widetilde{V}_{(h,\tau,\theta)}[\hat{r}(\mathrm{X})]\succeq J(\theta)[\widetilde{G}_n^{(h,\tau)}(\theta)]^{-1}J(\theta)^{\top},
   \end{equation}where $J(\theta)$ is the $k\times k $-Jacobian matrix of $r(\theta)$. Furthermore, equality holds in \eqref{eqn:t_CRLB} if and only if
   \begin{equation}\label{eqn:equality_condition_t_crlb}
     \hat{r}(\mathrm{x})-\widetilde{\mathbb{E}}_{(h,\tau,\theta)}[\hat{r}(\mathrm{X})]=J(\theta)\big[\widetilde{G}_n^{(h,\tau)}(\theta)\big]^{-1}\frac{1}{\widetilde{p}_{(h,\tau,\theta)}(\mathrm{x})}\frac{\partial p_{\theta}(\mathrm{x})}{\partial \theta}.
\end{equation}
\end{theorem}
\begin{proof}\
   The proof follows similarly to that of Theorem~\ref{thm:h_CRLB}.
\end{proof}
\begin{remark}
    For $h(t)=t\log t-t+1$ and $\tau(t)=t$,
the $\widetilde{(h,\tau)}$-CR inequality coincides
with the classical CR inequality \eqref{eqn:CRLB},
and the equality condition
\eqref{eqn:equality_condition_t_crlb} reduces
to the corresponding classical equality condition.
\end{remark}
\begin{remark}
    For $h(t)=\varphi(t)$, $\tau(t)=t$, the $\widetilde{(h,\tau)}$-CR inequality coincides
with the $\varphi$-CRLB established in \cite{dhadumia2026cramer}.
\end{remark}
 In the following, we establish a necessary and sufficient condition for equality in \eqref{eqn:t_CRLB} in terms of $(h,\tau)$-exponential family. To do that we assume the following additional assumptions:
 \begin{enumerate}
   \item[(C)]  The support of $p_{\theta}$ is independent of $\theta$. That is, $\mathfrak{X}_{\theta}=\mathfrak{X}$ for all $\theta$. \vspace{.2cm}
     \item[(D)] For each $\mathrm{x}\in \mathfrak{X}$, $(h'\circ\tau)\bigl(p_\theta(\mathrm{x})\bigr)$ is twice continuously differentiable in $\mathbb{\theta}$.\vspace{.2cm}
    \item[(E)] $\widetilde{V}_{(h,\tau,\theta)}[\hat{r}(\mathrm{X})]$ is positive definite for all $\theta$.
 \end{enumerate}
We prove the following theorem for the scalar case ($k=1$), although it holds for the vector parameter.
\begin{theorem}\label{thm:vector_h,tau_family_equality}
     Let $\mathrm{X}=(X_1,\ldots,X_n)$ be distributed according to $p_{\theta}$. Let $\hat{r}(\mathrm{X})$ be an unbiased estimator of $r(\theta)$ with respect to $p_{\theta}$. Suppose $(I)-(II)$ and $(A)-(E)$ are satisfied.  Then $\hat{r}(\mathrm{X})$ attains \eqref{eqn:t_CRLB} if and only if
\begin{equation}\label{eqn:vector_h,tau_family_equality}
    (h'\circ\tau)\bigl(p_\theta(\mathrm{x})\bigr) =
 C(\mathrm{x}) + Z(\theta)+\lambda(\theta)^{\top}\hat{r}(\mathrm{x})\qquad
\text{for }\mathrm{x}\in\mathcal{X},
  \end{equation} for some functions $C:\mathcal{X}\to \mathbb{R}$ and $Z\in C^1(\Theta)$, and $\lambda=[\lambda_1,\ldots,\lambda_k]^{\top}$, where $\lambda_i\in C^1(\Theta)$, $i=1,\ldots,k$, such that the $k\times k$ Jacobian matrix $R(\theta)\coloneqq\frac{\partial\lambda(\theta)}{\partial\theta}\coloneqq[\partial_j\lambda_i(\theta)]_{i,j=1}^{k}$ is invertible.
\end{theorem}
\begin{proof}\
    Assume that equality holds in \eqref{eqn:t_CRLB}. Then, from \eqref{eqn:equality_condition_t_crlb}, we have
\begin{equation}\label{eq:001}
    \hat{r}(\mathrm{x}) = a(\theta)\frac{1}{\widetilde{p}_{(h,\tau,\theta)}(\mathrm{x})}\frac{\partial p_{\theta}(\mathrm{x})}{\partial \theta}+b(\theta),
\end{equation}
with \begin{equation*}
    a(\theta) =
 r'(\theta)[\widetilde{G}_n^{(h,\tau)}(\theta)]^{-1} \quad\text{and}\quad b(\theta) = \widetilde{\mathbb{E}}_{(h,\tau,\theta)}[\hat{r}(\mathrm{X})].
\end{equation*}
Since $\widetilde{V}_{(h,\tau,\theta)}[\hat{r}(\mathrm{X})]>0$, it follows that
$\hat{r}(\mathrm{x})$ is not constant. Hence, $a(\theta)\neq 0$.
Therefore, from \eqref{eq:001}, we obtain
\begin{equation}\nonumber
   \frac{1}{\widetilde{p}_{(h,\tau,\theta)}(\mathrm{x})}\frac{\partial p_{\theta}(\mathrm{x})}{\partial \theta}
    =a(\theta)^{-1}\hat{r}(\mathrm{x})
    -a(\theta)^{-1}b(\theta).
\end{equation} This implies 
\begin{equation}\nonumber
   h''(\tau(p_{\theta}(\mathrm{x})))\tau'(p_{\theta}(\mathrm{x}))\frac{\partial p_{\theta}(\mathrm{x})}{\partial \theta}=k(\theta)a(\theta)^{-1}\hat{r}(\mathrm{x})-k(\theta)a(\theta)^{-1}b(\theta), 
\end{equation} where $k(\theta)\coloneqq1/\int [h''(\tau(p_{\theta}(\mathrm{y})))\tau'(p_{\theta}(\mathrm{y}))]^{-1}d\mu(\mathrm{y})$.
Integrating the above equation with respect to $\theta$ yields
\eqref{eqn:vector_h,tau_family_equality}.

Conversely, let us assume that \eqref{eqn:vector_h,tau_family_equality} holds. That is,
\begin{equation*}
 h'(\tau(p_{\theta}(\mathrm{x})))
= C(\mathrm{x}) + Z(\theta) + \lambda(\theta) \hat{r}(\mathrm{x}),\quad \mathrm{x}\in \mathcal{X}.
\end{equation*} Then
\begin{equation}\label{eqn:derivative_of_h_function}
    \frac{\partial}{\partial \theta} h'(\tau(p_{\theta}(\mathrm{x})))=Z'(\theta)+  \lambda'(\theta)\hat{r}(\mathrm{x}).
\end{equation}
Taking expectation on \eqref{eqn:derivative_of_h_function} with respect to $\widetilde{p}_{(h,\tau,\theta)}$, we get
\begin{equation}\nonumber
    \int  \frac{\partial}{\partial \theta} h'(\tau(p_{\theta}(\mathrm{x})))\widetilde{p}_{(h,\tau,\theta)}(\mathrm{x})d\mu(\mathrm{x})=Z'(\theta)+\lambda'(\theta)\widetilde{\mathbb{E}}_{(h,\tau,\theta)}[\hat{r}(\mathrm{X})].
\end{equation} That is,
\begin{equation}\nonumber
    k(\theta)\int \frac{\partial p_{\theta}(\mathrm{x})}{\partial\theta}d\mu(\mathrm{x})=Z'(\theta)+\lambda'(\theta)\widetilde{\mathbb{E}}_{(h,\tau,\theta)}[\hat{r}(\mathrm{X})].
\end{equation} This yields,
\begin{equation}
   Z'(\theta)=-\lambda'(\theta)\widetilde{\mathbb{E}}_{(h,\tau,\theta)}[\hat{r}(\mathrm{X})].
\end{equation} From \eqref{eqn:derivative_of_h_function}, we get 
\begin{equation}\nonumber
     \frac{\partial}{\partial \theta} h'(\tau(p_{\theta}(\mathrm{x})))=\lambda'(\theta)(\hat{r}(\mathrm{x})-\widetilde{\mathbb{E}}_{(h,\tau,\theta)}[\hat{r}(\mathrm{X})]).
\end{equation} Since $\lambda'(\theta)\neq0$, we get 
\begin{equation}\label{eqn:mx1}
    [\lambda'(\theta)]^{-1}  \frac{\partial}{\partial \theta} h'(\tau(p_{\theta}(\mathrm{x})))=\hat{r}(\mathrm{x})-\widetilde{\mathbb{E}}_{(h,\tau,\theta)}[\hat{r}(\mathrm{X})].
\end{equation}
Observe that
\begin{align*}
    \widetilde{G}_n^{(h,\tau)}(\theta)&=\int \frac{1}{\widetilde{p}_{(h,\tau,\theta)}(\mathrm{x})}\frac{\partial p_{\theta}(\mathrm{x})}{\partial \theta}\frac{\partial p_{\theta}(\mathrm{x})}{\partial \theta}d\mu(\mathrm{x})\\
    &=\frac{1}{k(\theta)}\int h''(\tau(p_{\theta}(\mathrm{x})))\tau'(p_{\theta}(\mathrm{x}))\frac{\partial p_{\theta}(\mathrm{x})}{\partial \theta}\frac{\partial p_{\theta}(\mathrm{x})}{\partial \theta}d\mu(\mathrm{x})\\
    &\overset{(a)}{=}\frac{1}{k(\theta)}\int (Z'(\theta)+\lambda'(\theta)\hat{r}(\mathrm{x}))\frac{\partial p_{\theta}(\mathrm{x})}{\partial \theta}d\mu(\mathrm{x})\\
    &=\frac{1}{k(\theta)}\lambda'(\theta)J(\theta),
\end{align*}where $(a)$ follows from \eqref{eqn:derivative_of_h_function}. That is,
\begin{equation}\label{eqn:mx2}
    [\lambda'(\theta)]^{-1}=\frac{1}{k(\theta)}J(\theta)[   \widetilde{G}_n^{(h,\tau)}(\theta)]^{-1}
\end{equation} By substituting \eqref{eqn:mx2} in \eqref{eqn:mx1}, we get \eqref{eqn:equality_condition_t_crlb}. This completes the proof.
\end{proof}
\begin{remark}
  For $h(t)=t\log t-t+1$ and $\tau(t)=t$,
\eqref{eqn:vector_h,tau_family_equality} reduces to
the exponential family, which characterizes
equality in the classical CR inequality
\cite{wijsman1973attainment,fend1959attainment,muller1989attainment}.
\end{remark}
\begin{remark}
   For $h(t)=\varphi(t)$ and $\tau(t)=t$, \eqref{eqn:vector_h,tau_family_equality} reduces to the $\varphi$-family, which characterizes equality in the $\varphi$-CRLB \cite[Theorem~6]{dhadumia2026cramer}.
\end{remark}
\begin{remark}
The proof of the $\widetilde{(h,\tau)}$-CR inequality in Theorem~\ref{thm:t_CRLB} is based on the Cauchy--Schwarz inequality.
An information-geometric proof can be obtained
by considering the following Bregman divergence:
\begin{equation}
\begin{aligned}\label{eqn:F_BRegman_div}
D_F(p_{\theta},p_{\theta'})
\coloneqq \int\Big[
&F(p_{\theta}(\mathrm{x}))
-F(p_{\theta'}(\mathrm{x}))-F'(p_{\theta'}(\mathrm{x}))
\big(p_{\theta}(\mathrm{x})-p_{\theta'}(\mathrm{x})\big)
\Big]\,d\mu(\mathrm{x}),
\end{aligned}
\end{equation}
where $F$ is an antiderivative of $h'\circ\tau$ on an open interval $I\subset \mathbb{R}$.
We omit this alternative proof, as it proceeds analogously
to the information-geometric proof of the $(h,\tau)$-CR
inequality using the $(h,\tau)$-divergence, with
$\widetilde{p}_{(h,\tau,\theta)}$ replacing
$p_{(h,\tau,\theta)}$.
\end{remark}
The following example illustrates that an estimator from an $(h,\tau)$-exponential family that yields strict inequality
in the $(h,\tau)$-CR inequality attains equality in the
$\widetilde{(h,\tau)}$-CR inequality.   
\begin{example}\label{example:suff_attain_t_CRLB}
Consider $p_{\theta}$ as defined in \eqref{eqn:Ber_alpha_form}
of Example~\ref{example:inequality_equality}, which forms an $(h,\tau)$-exponential family with $h(t)=\frac{2}{3}t^{3/2}$, $\tau(t)=t^2$
(see Remark~\ref{remark:alpha_family_ofBer}). 
In this case, $\widetilde{p}_{(h,\tau,\theta)}(\mathrm{x})=\frac14$
for all $\mathrm{x}$, and
\[
\widetilde{V}_{(h,\tau,\theta)}[\hat{r}(\mathrm{X})]
=
J(\theta)[\widetilde{G}_n^{(h,\tau)}(\theta)]^{-1}J(\theta)
=
\frac12.
\]
Hence, the estimator
$\hat r(\mathrm{X})=X_1+X_2-\frac{1}{2}$
attains the $\widetilde{(h,\tau)}$-CR inequality for $\theta$. In contrast, the \((h,\tau)\)-CR inequality is strict for the same estimator.
\end{example}

Note that from the Factorization Theorem \cite[Th.~6.2.6]{casella2024statistical}, in \eqref{eqn:Ber_alpha_form}, $\hat{r}(\mathrm{X})=X_1+X_2-\frac{1}{2}$ is sufficient for $\theta$. In Example~\ref{example:suff_attain_t_CRLB}, we have seen that $\hat{r}(\mathrm{X})$ attains the $\widetilde{(h,\tau)}$-CR inequality.
The following theorem extends this result to the general $(h,\tau)$-exponential family.
\begin{theorem}\label{thm:attainment_sufficient_h_tau_family}
     Suppose $p_{\theta}$ forms an $(h,\tau)$-exponential family as in \eqref{eqn:vector_h,tau_family_equality} with $C(\mathrm{x})=C\in\mathbb{R}$ for all $\mathrm{x}\in\mathcal{X}$. Then \(\hat{r}(\mathrm{X})\) is a sufficient statistic for \(\theta\). Moreover, if $r(\theta)=\mathbb{E}_{\theta}[\hat{r}(\mathrm{X})],$ then \(\hat{r}(\mathrm{X})\) attains the $\widetilde{(h,\tau)}$-CR inequality for the parameter \(r(\theta)\). 
\end{theorem}
\begin{proof}\
    Follows from the Factorization Theorem and Theorem~\ref{thm:vector_h,tau_family_equality}.
\end{proof}
\begin{remark}
Theorem~\ref{thm:attainment_sufficient_h_tau_family}
extends the classical result that a sufficient statistic
of an exponential family attains the CR inequality
for its expectation parameter to $(h,\tau)$-exponential family  \cite[Eq.~(3.13)]{crrao1945}.
\end{remark} We present another example in which equality in the $\widetilde{(h,\tau)}$-CR inequality is attained by the sufficient statistic of an $(h,\tau)$-exponential family.
\begin{example}
    Let $\mathrm{X}=(X_1,\ldots,X_n)$ be distributed according to the multivariate Student $t$-distribution $t_{n,\nu}(\mu\mathbf{1}_n,\sigma^2\mathbb{I}_n)$, $\mu\in\mathbb{R}$, $\sigma^2\in (0,\infty)$, with $\nu>2$, where
    \begin{equation}
        \label{eqn:mult_t_dist}
       p_{\mu,\sigma^2}(\mathrm{x})= {\displaystyle {N_{\nu,\sigma^2}}\left[1+\frac{1}{\nu\sigma^2}\sum_{i=1}^n(x_i-\mu)^2\right]^{-(\nu +n)/2}},\quad \mathrm{x}=[x_1,\ldots,x_n]^{\top}\in\mathbb{R}^n,
    \end{equation} where
$$N_{\nu,\sigma^2}=\frac {\Gamma \left[(\nu +n)/2\right]}{\Gamma (\nu /2)\nu ^{n/2} \pi^{n/2}{\sigma}^{n}}.$$
\eqref{eqn:mult_t_dist} forms an $(h,\tau)$-exponential family with $h(t)=-\frac{\nu+n}{\nu+n-2}t^{\frac{\nu+n-2}{\nu+n}}$ and $\tau(t)=t$, $C(\mathrm{x})=0$, $Z(\theta)=-N_{\nu,\sigma^2}^{-\frac{2}{\nu+n}}
-\frac{n\mu^2}{\nu\sigma^2} N_{\nu,\sigma^2}^{-\frac{2}{\nu+n}},$ $\lambda(\theta) =\left[
\frac{2n\mu}{\nu\sigma^2}N_{\nu,\sigma^2}^{-\frac{2}{\nu+n}},-
\frac{n}{\nu\sigma^2} N_{\nu,\sigma^2}^{-\frac{2}{\nu+n}}
\right]^\top$, $T(\mathrm{x})=
\left[
\frac{1}{n}\sum_{i=1}^n x_i,\,
\frac{1}{n}\sum_{i=1}^n x_i^2
\right]^{\top}$.
It can be shown that $T(\mathrm{x})$ is a sufficient statistic for $[\mu,\sigma^2]^{\top}$ and $\mathbb{E}_{p_{(\mu,\sigma^2)}}[T(\mathrm{X})]=[\mu,\mu^2+\frac{\nu}{\nu-2}\sigma^2]^{\top}$. By Theorem~\ref{thm:attainment_sufficient_h_tau_family}, $T(\mathrm{X})$ attains the $\widetilde{(h,\tau)}$-CR inequality for its expected parameter.
\end{example}
\begin{remark}
The $(h,\tau)$-CR inequality and $\widetilde{(h,\tau)}$-CR inequality share structural similarities with the generalized Cramér–Rao bound introduced by Jan Naudts \cite{naudts2004estimators}, particularly through the appearance of escort distribution in their formulation. While Naudts considers a general escort distribution $P_{\theta}$, the specific escort distributions
$p_{(h,\tau,\theta)}$ and
$\widetilde{p}_{(h,\tau,\theta)}$ emerge naturally
from our framework. Although both inequalities fall within the general escort-based framework of Jan Naudts, our motivation and approach to deriving these inequalities differ from those of Naudts.
\end{remark}
\section{Summary }\label{sec:5}
In the literature, the $(h,\tau)$-divergence and the $(h,\tau)$-exponential family have been established. Upon choosing $(h(t),\tau(t))=\bigl(t\log t-t+1,\,t\bigr)$, they reduce to the KL divergence and the classical exponential family, respectively. Amari and Nagaoka gave an information-geometric proof of the classical CR inequality based on the Fisher information metric, which can be obtained from the KL divergence using Eguchi's theory. In this paper, we extended this theory to the $(h,\tau)$-divergence and established the $(h,\tau)$-CR inequality. To this end, we first applied Eguchi's theory to the $(h,\tau)$-divergence to obtain the corresponding Riemannian metric and then apply Amari–Nagaoka's theory based on this metric. We also gave a proof of this CR inequality using Cauchy--Schwarz inequality, enabling the application of this bound to more general support sets.

The classical result that equality holds in CR inequality under exponential family does not generally extend to $(h,\tau)$-CR inequality under $(h,\tau)$-exponential family. Consequently, we established a CR inequality that is attained under $(h,\tau)$-exponential family. This result extends the classical attainment of the CR inequality by the sufficient statistic of the exponential family to the $(h,\tau)$-exponential family. Finally, we identified a Bregman divergence \eqref{eqn:F_BRegman_div} associated with the $\widetilde{(h,\tau)}$-CR inequality, providing an alternative information-geometric approach to its derivation.
\section{Acknowledgements}
Satyajit Dhadumia is supported by a Council of Scientific \& Industrial Research (CSIR) fellowship of the Department of Scientific \& Industrial Research, Ministry of Science and Technology, Government of India.
\bibliographystyle{plain}
    \bibliography{ref}

@article{crrao1945,
  title={Information and the accuracy attainable in the estimation of statistical parameters},
  author={Rao, C. Radhakrishna},
  journal={Bull. Calcutta Math. Soc},
  volume={37},
  number={3},
  pages={81--91},
  year={1945}
}

@article{Eguchi1992,
    AUTHOR = "Eguchi, Shinto",
     TITLE = "Geometry of minimum contrast",
   JOURNAL = "Hiroshima Math. J.",
  FJOURNAL = "Hiroshima Mathematical Journal",
    VOLUME = "22",
      YEAR = "1992",
    NUMBER = "3",
     PAGES = "631--647",
      ISSN = "0018-2079",
   MRCLASS = "62B10 (62A99)",
}

@book{AmariHiroshi2000,
    AUTHOR = "Amari, S and Nagaoka, Hiroshi",
     TITLE = "Methods of information geometry",
    SERIES = "Translations of Mathematical Monographs",
    VOLUME = "191",
      NOTE = "Translated from the 1993 Japanese original by Daishi Harada",
 PUBLISHER = "American Mathematical Society, Providence, RI; Oxford
              University Press, Oxford",
      YEAR = "2000",
     PAGES = "x+206",
      ISBN = "0-8218-0531-2",
   MRCLASS = "62F12 (53B99 81Q70 94A15)",
}

@article{naudts2004estimators,
  title={Estimators, escort probabilities, and exponential families in statistical physics.},
  author={Naudts, Jan},
  journal={JIPAM. Journal of Inequalities in Pure \& Applied Mathematics},
  volume={5},
  number={4},
  pages={Paper--No},
  year={2004},
  publisher={Victoria University, School of Communications and Informatics}
}

@book{casella2024statistical,
  title={Statistical inference},
  author={Casella, George and Berger, Roger},
  year={2024},
  publisher={Chapman and Hall/CRC}
}

@article{Roy2019DensityPower,
  title   = {Density Power Downweighting and Robust Inference: Some New Strategies},
  author  = {Roy, Saptarshi and Chakraborty, Kaustav and Bhadra, Somnath and Basu, Ayanendranath},
  journal = {Journal of Mathematics and Statistics},
  volume  = {15},
  pages   = {333--353},
  year    = {2019},
  doi     = {10.3844/jmssp.2019.333.353}
}

@inproceedings{dhadumia2025generalized,
  title={A Generalized Cram{\'e}R-Rao Bound Using Information Geometry},
  author={Dhadumia, Satyajit and Kumar, M Ashok},
  booktitle={2025 IEEE International Symposium on Information Theory (ISIT)},
  pages={1--6},
  year={2025},
  organization={IEEE}
}

@article{fend1959attainment,
  title={On the attainment of Cramer-Rao and Bhattacharyya bounds for the variance of an estimate},
  author={Fend, AV},
  journal={The Annals of Mathematical Statistics},
  pages={381--388},
  year={1959},
  publisher={JSTOR}
}

@article{wijsman1973attainment,
  title={On the attainment of the Cram{\'e}r-Rao lower bound},
  author={Wijsman, Robert A},
  journal={The Annals of Statistics},
  pages={538--542},
  year={1973},
  publisher={JSTOR}
}

@article{muller1989attainment,
  title={On the attainment of the Cram{\'e}r-Rao bound in {L}r-differentiable families of distributions},
  author={M{\"u}ller-Funk, Ulrich and Pukelsheim, Friedrich and Witting, Hermann},
  journal={The Annals of Statistics},
  pages={1742--1748},
  year={1989},
  publisher={JSTOR}
}

@article{dhadumia2026cramer,
  title={Cram{\'e}r--{R}ao Type Bounds for Generalized Exponential Families},
  author={Dhadumia, Satyajit and Kumar, M. Ashok},
  journal={arXiv preprint arXiv:2610.27679},
  year={2026}
}

@article{MATSUZOE2027131088,
title = {Divergence and deformed exponential family},
journal = {Journal of Mathematical Analysis and Applications},
volume = {567},
number = {1},
pages = {131088},
year = {2027},
issn = {0022-247X},
doi = {https://doi.org/10.1016/j.jmaa.2026.131088},
url = {https://www.sciencedirect.com/science/article/pii/S0022247X26007006},
author = {Hiroshi Matsuzoe and Asuka Takatsu}
}

\end{document}